\documentclass[]{siamart220329}
\usepackage{amsmath}
\usepackage{amssymb}
\usepackage{graphicx}
\usepackage{calc}

\usepackage{pgfplots}
\usepgfplotslibrary{external}
\usepackage{hyperref}
\usepackage{subcaption}

\newcommand{\R}{\mathbb{R}}

\newtheorem{remark}{Remark}
\newtheorem{example}{Example}

\begin{document}
	
	\title{A hierarchy of kinetic discrete-velocity  models for traffic flow derived from a non-local Prigogine-Herman model}
	
	\author{R. Borsche\footnotemark[1] 
		\and  A. Klar\footnotemark[1] \footnotemark[2]}
	\footnotetext[1]{Technische Universit\"at Kaiserslautern, Department of Mathematics, Erwin-Schr\"odinger-Stra{\ss}e, 67663 Kaiserslautern, Germany 
		(\{borsche, klar\}@mathematik.uni-kl.de)}
	\footnotetext[2]{Fraunhofer ITWM, Fraunhoferplatz 1, 67663 Kaiserslautern, Germany} 
	
	\date{}
	

	\maketitle
	
	\begin{abstract}
		Starting from a non-local version of the Prigogine-Herman traffic model, we derive  a natural hierarchy of kinetic discrete velocity models   for traffic flow consisting of systems of quasi-linear hyperbolic equations with relaxation terms. 
		The hyperbolic main part of these models turns out to have several favourable features. In particular,  we determine Riemann invariants and  prove 
		richness and total linear degeneracy of the hyperbolic systems. 
		Moreover, a physically reasonable invariant domain is  obtained for all equations of the hierarchy.
		Additionally, we investigate the full relaxation system with respect to stability and persistence of periodic 
		(stop and go type) solutions
		and derive a condition for the appearance of such solutions.
		Finally, numerical results for various situations are presented, illustrating the analytical findings.
	\end{abstract}
	
	

	
	\section{Introduction and Motivation}
	
	Starting with the first macroscopic traffic model in \cite{Whi74}, there have been many approaches to a continuous modeling of traffic flow problems.
	Macroscopic models are usually based on scalar hyperbolic equations like 
	the above cited model or on systems of hyperbolic equations like in  \cite{pay79,AR,Zhang}.
	For discussions and extensions of these models see, for example, \cite{AKMR,Deg,Ber,G,R}.
	On the other hand, kinetic equations have also been widely used as a tool to model traffic flow problems. Starting with the work in \cite{PH71}  different extensions and amplifications can be found in  \cite{Hel95B,KW981,Nel95,PF75} or \cite{KW97,KW00}.   
	For more recent works, we refer  to \cite{FT13} for a fully discrete kinetic approach, to \cite{HPRV} for a BGK-type traffic model, to \cite{CT,DT,DTZ,BK22} for relations between kinetic and macroscopic models and to \cite{BD11} for a mathematically oriented review and further references.
	
	A naive use of classical  kinetic  equations  in the case of traffic flow leads to well known problems:
	physically given invariant domains (a state space restricted to positive and bounded velocities and densities) are not respected and information transport against the flow direction
	is not possible.
	The second issue has been treated in several works  introducing 
	some additional mechanism   into the equations.
	The Enskog type  kinetic approaches in \cite{KW97} use a non-locality to deal with these problem.
	Similarly, the discrete model developed in \cite{FT13} 
	includes the non-locality via a limiter term into the free flow part of the kinetic equation. 	
	The model in \cite{HPRV} uses a BGK-type approximation and a suitable modification to allow for a proper information transport.

	In this work we  derive  a hierarchy of kinetic discrete velocity traffic models from a non-local version of the original Priogine-Herman kinetic traffic model \cite{PH71}.
	In fluid- or gas-dynamics  kinetic discrete velocity models have been investigated in many works, see
	\cite{Illner} for a review. 
	Consider as an example a  classical continuous BGK-type kinetic equation for the kinetic distribution function $f=f(x,v,t) \ge 0$
	\begin{align*}
		\begin{aligned}
			\partial_t f + v \partial_x f = -\frac{1}{\epsilon} \left(f-f^e(\rho,v) \right)
		\end{aligned}
	\end{align*}
	with $ \int f(v) d v = \rho = \int f^e(\rho,v ) d v $.
	An associated discrete velocity model has the distribution functions $f_i,i=0, \cdots N$ for velocities $ v_0  < v_1  < \cdots < v_i < \cdots <  v_N $.
	Let the discrete equilibrium functions be denoted by $ f^e_i (\rho) $,  the density $\rho$ be given by $\rho = f_0 + \cdots + f_N =f^e_0 + \cdots   +f^e_N$ and the mean flux by $q = v_0 f_0 + v_1 f_1 +\cdots +  v_Nf_N$. 
	Then, with the equilibrium flux $F(\rho) = v_0 f^e_0 + v_1 f^e_1 + \cdots + v_N f^e_N $, the discrete-velocity model is given by 
	\begin{align}
		\label{dvm0}
		\begin{aligned}
			\partial_t f_i + v_i \partial_x f_i = -\frac{1}{\epsilon} \left(f_i-f^e_i(\rho) \right).
		\end{aligned}
	\end{align}
	The associated conservation law for the density for small values of $\epsilon$ is formally given by 
	$$
	\partial_t \rho  + \partial_x F(\rho) =0.
	$$ 
	Applying this framework naively to traffic flow modeling, there are two simple observations related to each other.   
	First, since in the traffic situation all velocities in the kinetic model are positive, i.e. $v_i \ge 0,  \; i=0, \ldots,N$, there is no backward information transport, as, for example,
	in the case of a backward traveling traffic jam, since there are no negative wave speeds.
	The fact that there cannot be negative wave speeds also prevents a possible convergence to the limiting conservation law
	$
	\partial_t \rho + \partial_x F(\rho) =0,
	$
	if there are  situations with density $\rho$, such that  $F^\prime (\rho )$ is  negative, due to the sub characteristic condition,  see e.g. \cite{Liu87}. 
	
	Second,  since the  density $\rho$ in a discrete velocity traffic model is limited by  a maximal (bumper-to-bumper) density, the physically reasonable region for such a model is  given by the N-simplex $0 \le f_i$ and $\rho =f_0+  \cdots + f_N\le \rho_{max}$. 
	However, this is not an invariant domain for the above equation \eqref{dvm0}. 
	This can be easily seen looking at the hyperbolic part and considering a suitable Riemann problem.
	The bound on $\rho $ is not satisfied by the solutions of the equation and one obtains solutions with
	densities exceeding the bumper-to-bumper density. See  the related discussion in \cite{Deg} for macroscopic traffic models.

	In the present paper we introduce a natural new 
	class of discrete-velocity models for traffic flow derived from a non-local version of  the classical Prigogine-Herman model
	and allowing a more rigorous investigation than previous approaches.
	To derive the equations we proceed similarly as in \cite{BK18}, where a strongly simplified relaxation model has been derived and investigated.
	It will turn out that the resulting model is a rich and totally linear degenerate system of hyperbolic equations with relaxation.
	Such systems have many favourable properties and have been studied in detail in the literature, see, for example, \cite{serre,peng,Li,CHN09}.
	Moreover, this class of discrete velocity models resolves the above issues:
	it is shown  that the 
	invariant region is given by 
	the N-simplex discussed above and that the equations allow for negative wave speeds. Moreover,
	the full relaxation model exhibits  stable and unstable flow situations allowing the modeling of stable traffic and 
	of stop-and-go type instabilities in traffic flow.

	Our starting point is  a classical continuous kinetic traffic flow equation or, more exactly,  a modified version of the 
	spatially inhomogeneous original Prigogine-Herman equation.
	We proceed with a discretization in velocity space of these equations. A   detailed investigation of the resulting hyperbolic  parts of the equations of the discrete hierarchy  gives  a rich hyperbolic system  and the invariant domain discussed above.
	Moreover, the equations together with the relaxation term are considered and a  condition  for stability respectively instability and appearance of persistent periodic (stop and go) waves is given following the approach in  \cite{G04,GKR03}. See also  \cite{Flynn} for a discussion of travelling wave solutions for macroscopic traffic equations.

	The paper is organized in the following way. 
	In section \ref{kineticmodel} we discuss a continuous kinetic traffic model based on the classical Prigogine Herman equation.
	In section \ref{discretevelocity}  we derive a discrete velocity kinetic model  from the  continuous kinetic traffic equation.
	In the subsequent sections \ref{discreteinvest1} and \ref{discreteinvest2} the hyperbolic parts of the discrete kinetic  equations are investigated in detail.
	In particular, we prove that the system is a rich and  totally linear degenerate hyperbolic system.
	Section \ref{relax} discusses the full relaxation model and the above mentioned stability issues for continuous and discrete equations.
	Finally,
	numerical results are presented in Section \ref{numericalresults}.

	\section{A continuous kinetic model for traffic flow}
	\label{kineticmodel}
	
	Our starting point is a modified version of the classical kinetic Prigogine-Herman equation for traffic flow, see  \cite{KW97,KW981,KW00,Hel95B,Nel95,PF75,PH71}. We proceed similarly as in \cite{BK18}.

	\subsection{Derivation of the continuous kinetic model}
	The kinetic model is based on an underlying microscopic model with breaking interactions, where the driver at $x$ with velocity $v$ reacts to his predecessor at $x+H$ with velocity $\hat v < v$. Here   $H> 0$ is  a minimal distance between the vehicles.
	The new velocity resulting of this braking interaction is equal to  the velocity of the leading car.
	Then, $\frac{1}{H}$ is the maximal density of vehicles on the road and the interaction strength is modulated by a factor $\frac{1}{1-H \rho}$ which can be  derived from the two particle correlation function of the vehicles, compare \cite{KW00}. 
	One obtains for 
	$t \in \R^+$, $x \in \R$, $v \in[0,v_{max}]$ the following equation for the distribution function $f=f(x,v,t)\ge 0$ 
	\begin{align*}
		\partial_t f + v \partial_x f = J (f) =  J_B(f) + J_A(f),
	\end{align*}
	where $J_B$ describes the interactions with respect to braking and $J_A$ additional interactions, like those related to acceleration.
	The braking term $J_B (f)$ is given as in  \cite{PH71,BKK} by 
	\begin{align*}
		J_B (f) = \frac{1}{1-H \rho }\int_{\hat v > v}  (\hat v -v)  f(x,\hat v) f(x+H, v ) d \hat v\\
		- \frac{1}{1-H\rho }\int_{\hat v < v } ( v - \hat v)  f(x, v) f(x+H, \hat v ) d \hat v\ 
	\end{align*}
	with $\rho= \int_0^{v_{max}} f(v) d v \le \rho_{max} = \frac{1}{H}$.

	Nondimensionalization using the new variables $\tilde x = \frac{x }{H }, \tilde t = \frac{v_{max}}{H } t$, $\tilde v =\frac{v}{v_{max}}$ and the new distribution function
	$\tilde f(\tilde x, \tilde v, \tilde t) = \frac{v_{max}}{\rho_{max}} f(x,v,t)= H v_{max} f (x,v,t)$ gives, if we neglect the tilde-notation,
	the equation
	\begin{align}
		\label{kineticunscaled}
		\partial_t f + v \partial_x f =   J_B(f)  + J_A (f) 
	\end{align}
	with
	\begin{align*}
		J_B (f) = \frac{1}{1-\rho }\int_{\hat v > v}  (\hat v -v)  f(x,\hat v) f(x+1, v ) d \hat v\\
		- \frac{1}{1-\rho }\int_{\hat v < v } ( v - \hat v)  f(x, v) f(x+1, \hat v ) d \hat v\ ,
	\end{align*}
	where   $t \in \R^+$, $x \in \R$, $v \in[0,1]$  and $\rho= \int_0^1 f(v) d v  \le 1$.

Now, we split the braking term into two terms containing  the local and non-local braking effects, respectively.
The local braking term is 
\begin{align*}
	J_B^L (f) &= \frac{1}{1-\rho }\int_{\hat v > v}  (\hat v -v)  f(x,\hat v) f(x, v ) d \hat v\\
	&- \frac{1}{1-\rho }\int_{\hat v < v } ( v - \hat v)  f(x, v) f(x, \hat v ) d \hat v \\
	&=\frac{1}{1-\rho }	\int (\hat v-v) f(v) f(\hat v)d \hat v 
	=  f(v) \frac{q-v \rho}{1-\rho }
\end{align*}
with $q= \int_0^1 v f(v) d v $.
\begin{remark} The term  $J_B^L$ is the classical breaking term for kinetic problems as  in the work of Prigogine and Herman
	\cite{PH71}.
\end{remark}

The  term containing the  non-local effects due to braking interactions is then denoted by $J_B^{NL}(f)$. Thus,
$$
J_B^{NL} (f) =  J_B (f) - J_B^L (f) .
$$
The total local effects are then described by the  combined term
\begin{align*}
	J^L(f) =  J_A (f) +J_B^L (f) .
\end{align*}
To simplify the situation we approximate this term by a standard relaxation term
\begin{align*}
	J^L(f) \approx J_R (f) =  - \frac{1}{T} \left( f - f^e (\rho,v) \right)\ 
\end{align*}
with an equilibrium function  $f^e$ fulfilling $\int f^e (\rho,v) dv = \rho = \int f (v) dv$. The first moment  of $f^e$   is  $\int v f^e(\rho,v) dv = F(\rho)$, where 
$F(\rho) , 0 \le \rho \le 1$ is called the fundamental diagram.
$f^e $ is chosen in such a way, that  $V(\rho)= \frac{F(\rho)}{\rho}$ is monotone decaying with $V(0) = 1 $  and 
$V(1) =0$. 
This yields functions $F$ fulfilling  $F(0) = F(1) =0$ and $F^{\prime}(\rho) \le 1$. 
Moreover, in the following, we use $E=E(\rho) $ to denote the second moment of $f^e$, i.e. $E =\int v^2 f^e(\rho,v)dv $.

\begin{example}
	\label{excont}
	As an example for the equilibrium function $f^e$, which  will also be used for numerical computations, we consider  a class of equilibrium functions $f^e$ depending only on the first 3 moments $\rho,F,E$, i.e.
	\begin{align*}
		f^e(\rho,v) = \delta_0 (v)  \left(\rho-F - (\gamma-1) \frac{F-E}{1-\lambda} \right) +  \frac{F-E}{1-\lambda} \frac{g(v)}{v}+ \delta_1(v) \left(  F-\frac{F-E}{1-\lambda}\right)
	\end{align*}
	with a smooth function $g \ge 0$, $\int_0^1  g(v) dv =1$. Here $\lambda = \int_0^1 v g(v) dv$ and 
	we require  $\gamma = \int_0^1 \frac{g(v)}{v} d v < \infty$.
	Note that realizability, i.e. positivity of $f^e$, requires
	\begin{align}
		\label{realcont}
		F- \frac{1-\lambda}{\gamma-1}  (\rho-F) \le E \; \mbox{and} \; \lambda F \le E \le F
	\end{align} 
	for all $0 < \rho < 1$.
	
\end{example}

To continue, the kinetic equation is now scaled with a hyperbolic space time scaling $t \rightarrow \frac{t}{\epsilon}$
and $x \rightarrow \frac{x}{\epsilon}$. 
This leads to 
\begin{align*}
	\partial_t f + v \partial_x f = \frac{1}{\epsilon }J_R (f) +
	\frac{1}{\epsilon } J_B^{NL} (f)
\end{align*}
with
\begin{align*}
	J_B^{NL} (f) &= \frac{1}{1-\rho }\int_{\hat v > v}  (\hat v -v)  f(x,\hat v)[f(x+\epsilon , v ) -f(x , v )] d \hat v\\
	&- \frac{1}{1-\rho }\int_{\hat v < v } ( v - \hat v)  f(x, v) [f(x+\epsilon , \hat v ) -f(x , \hat v )] d \hat v.
\end{align*}

\begin{remark}
	Note that in the fluid dynamic context, this is the classical hydrodynamic Enskog-Boltzmann scaling, see, for example, \cite{lacho}.
\end{remark}

Using a Taylor approximation we approximate the non-local term to order $\mathcal{O} (\epsilon^2)$ by the term $\epsilon J(f, \partial_x f ) $
with
\begin{align*}
	J(f, \partial_x f ) = \frac{ 1 }{1-\rho}
	\int_{\hat v > v}  ( \hat v -v )  f(\hat v)   \partial_x f (v)d \hat v 
	- \frac{1}{1-\rho } \int_{\hat v < v}  ( v-\hat v  ) f( v)   \partial_x f (\hat v)  d \hat v\  .
\end{align*}
Altogether,  we   obtain  the  kinetic problem
\begin{align}
	\label{kinbasic}
	\partial_t f + v \partial_x f - J (f, \partial_x f ) = \frac{1}{\epsilon }J_R (f) ,
\end{align}
which will be the starting point of our investigations.

\section{Discrete velocity  models for traffic flow}
\label{discretevelocity}

To obtain a hierarchy of discrete velocity models we discretize equation (\ref{kinbasic}) in velocity space with the (not necessarily equidistant) velocities $0 = v_0 <  v_1<v_2 < \cdots < v_{N-1} < v_N = 1$.
Note that the choices $v_0=0$ and $v_N = 1$ are necessary to span the full range of all possible velocities.
The associated discrete distribution functions we denote $f_i= f_i(x,t) \in [0,1], i=0, \ldots, N$ and define the discrete density and momentum by
\begin{align*}
	\rho = \sum_{i=0}^N f_i \qquad \mbox{and} \qquad
	q = \sum_{i=0}^N v_i f_i  .
\end{align*} 
For the discretization of the  equilibrium distribution we use  discrete 
values $f^e_i (\rho)$,  $i=0, \ldots, N$  fulfilling
$
\sum_{i=0}^N f^e_i (\rho) = \rho
$
and denote again
$
\sum_{i=0}^N v_i f^e_i (\rho) = F(\rho).
$

We obtain  the following discretization of the term $(1-\rho) J(f,\partial_x f)$.
\begin{align*}
	i=0:  & \; \sum_{j=1}^{N} v_j f_j \partial_x f_0\\
	i=1, \ldots,N-1:  &  \; \sum_{j=i+1}^{N} (v_j-v_i) f_j  \partial_x f_i - \sum_{j=0}^{i-1} (v_i-v_j) f_i \partial_x f_j  \\
	i=N:& \; - \sum_{j=0}^{N-1} (1-v_j) f_N  \partial_x f_j    .
\end{align*}
Altogether, this gives the hyperbolic  system
\begin{align}
	\label{discretekin}
	\begin{aligned}
		\partial_t f_0 + v_0 \partial_x f_0 - \frac{1}{1-\rho} \left[ \sum_{j=1}^{N} (v_j-v_0) f_j\partial_x f_0 \right] =J_R^0\\
		\partial_t f_i + v_i \partial_x f_i - \frac{1}{1-\rho} \left[\sum_{j=i+1}^{N} (v_j-v_i) f_j  \partial_x f_i - \sum_{j=0}^{i-1} (v_i-v_j) f_i \partial_x f_j  \right] =J_R^i\\
		\partial_t f_N + v_N \partial_x f_N - \frac{1}{1-\rho} \left[ - \sum_{j=0}^{N-1} (v_N-v_j) f_N  \partial_x f_j  \right]  =J_R^N
	\end{aligned}
\end{align}
for $i=1, \ldots,N-1$.
Here
\begin{align*}
	J_R^i = - \frac{1}{T} (f_i-f^e_i(\rho)) 
\end{align*}
for $i=0, \ldots,N$.
This is a quasilinear hyperbolic system 
$$
\partial_t f + A(f) \partial_x f= J_R(f),
$$
where the lower triangular matrix has the values
\begin{align*} 
	A_{ii} (f)= v_i  - \frac{1}{1-\rho} \left(\sum_{j=i+1}^{N} (v_j-v_i) f_j \right) 
	\ , \quad 
	A_{ij} (f)= \frac{1}{1-\rho}  (v_i - v_j) f_i  \ge 0
	\text{ for }j<i
	.
\end{align*}
The eigenvalues of $A$ are 
\begin{align*}
	\lambda_0&= - \frac{1}{1-\rho} \sum_{j=1}^{N} v_j f_j  \le  0,&&\lambda_N = 1, \\
	\lambda_i &= v_i- \frac{1}{1-\rho} \sum_{j=i+1}^{N} (v_j-v_i) f_j , &&i=1, \ldots, N-1
\end{align*}

\begin{example}
	\label{exN1}
	$N=1$.
	For  a two velocity model with the velocities $0 = v_0,  v_1 = 1$ 
	we obtain
	\begin{align*}
		\partial_t f_0 - \frac{1}{1-\rho }    f_1\partial_x f_0&= -\frac{1}{\epsilon} \left(f_0-\rho+F(\rho) \right)  \\
		\partial_t f_1 +  \partial_x f_1 + \frac{1}{1-\rho }    f_1\partial_x f_0&=  -\frac{1}{\epsilon} \left(f_1-F(\rho)\right).
	\end{align*}
	Note that this is the simplified relaxation model  derived in \cite{BK18}.
\end{example}

\begin{example}
	$N=2$.
	\label{exN2}
	For  the three  velocity model we obtain  with $\rho = f_0+f_1+f_2$ and $q=v_1f_1 +f_2$ and $e= v_1^2 f_1 + f_2$ and $E (\rho) = v_1^2 f^e_1 + f^e_2$, the following expressions
	$f_0 = \rho - f_1-f_2$, $f_1 = \frac{q-e}{v_1(1-v_1)}$, $f_2=\frac{e-v_1q}{1-v_1}$ and 
	\begin{align*}
		\partial_t f_0 - \frac{q}{1-\rho }    \partial_x f_0&= -\frac{1}{\epsilon} \left(f_0-\rho+ \frac{(1-v_1^2) F(\rho) - E(\rho) (1-v_1))}{v_1(1-v_1)}\right)  \\
		\partial_t f_1 + \left(v_1-   \frac{1-v_1}{1-\rho }   f_2 \right) \partial_x f_1 &+ \frac{v_1}{1-\rho }   f_1\partial_x f_0  =  -\frac{1}{\epsilon} \left(f_1-\frac{F(\rho)-E(\rho)}{v_1(1-v_1)}\right)\\
		\partial_t f_2 +  \partial_x f_2 + \frac{1-v_1}{1-\rho }    f_2\partial_x f_1 &+ \frac{1}{1-\rho }f_2 \partial_x f_0=  -\frac{1}{\epsilon} \left(f_2-\frac{v_1E(\rho)-v_1^2 F(\rho)}{v_1(1-v_1)})\right).
	\end{align*}
	Suitable functions for $F$ and $E$  lie in the realizability region
	\begin{align}
		\label{real}
		F- v_1(\rho-F) \le E,  \; v_1 F\le E \le F.
	\end{align}
\end{example}

\begin{example}
	\label{exN3}
	General $N$.	
	We obtain an example for general $N$ via  considering a discrete version of the equilibrium distribution in Example \ref{excont}.  
	That means we consider equilibrium functions of the following form.  
	Fix a function $F=F(\rho)$ and $E=E(\rho)$. Moreover, let $0 \le \alpha_i \le  1, i=1, \ldots , N-1$ and $\sum_{i=1}^{N-1} \alpha_i =1$. Moreover, we define $\lambda = \sum_{i=1}^{N-1} \alpha_i  v_i$ and $\gamma = \sum_{i=1}^{N-1} \frac{\alpha_i}{  v_i}$.  Then
	\begin{align*} 
		f_i^e(\rho,v) &= \frac{\alpha_i  }{v_i}\frac{ F-E}{1-\lambda}, &&i=1, \ldots, N-1,\\
		f_N^e(\rho,v) &= F- \sum_{i=1}^{N-1} v_i f_i^e\ , &&
		f_0^e(\rho,v) = \rho- \sum_{i=1}^{N}  f_i^e
	\end{align*}	
	Realizability requires again $F- \frac{1-\lambda}{\gamma-1}  (\rho-F) \le E$  and $\lambda F \le E \le F$ for all $0 < \rho < 1$ as in the continuous case.
	Note that for $N=2$ this is Example \ref{exN2}.
	
\end{example}

\section{Properties of the hyperbolic model}
\label{discreteinvest1}
In this and the following section we will give a thorough discussion of the hyperbolic part of  \eqref{discretekin} and show that the hyperbolic systems are, for all values of $N$, rich and totally linear degenerate, which leads to global existence of solutions, see \cite{peng,serre}.
We start with investigating the properties of the eigenvalues and eigenvectors of the system matrix.
Obviously  $\nabla \lambda_N = 0 $ and for  $i=0, \ldots, N-1$ we have with $\alpha = \frac{1}{1-\rho}$
\begin{align*}
	\nabla\lambda_i = \left(\begin{matrix}- \frac{1}{(1-\rho)^2} \sum_{j=i+1}^{N} (v_j-v_i) f_j 
		\\\vdots\\  
		- \frac{1}{(1-\rho)^2} \sum_{j=i+1}^{N} (v_j-v_i) f_j 
		\\
		- \frac{1}{1-\rho}  (v_{i+1}-v_{i} ) - \frac{1}{(1-\rho)^2} \sum_{j=i+1}^{N} (v_j-v_i) f_j    
		\\ \vdots \\ - \frac{1}{1-\rho}  (v_k -v_{i} ) - \frac{1}{(1-\rho)^2} \sum_{j=i+1}^{N} (v_j-v_i) f_j     \\ \vdots \\- \frac{H}{1-\rho}  (v_N-v_i)- \frac{1}{(1-\rho)^2} \sum_{j=i+1}^{N} (v_j-v_i) f_j  \end{matrix}\right)
	=  \alpha \left(\begin{matrix}\lambda_i-v_i
		\\\vdots\\  
		\lambda_i-v_i
		\\
		\lambda_i-v_{i+1}
		\\ \vdots \\  \lambda_i-v_{k}\\ \vdots \\ \lambda_i-v_{N}\end{matrix}\right).
\end{align*}
The eigenvectors  $r_i$ for the eigenvalues $\lambda_i, , i=0, \ldots, N$ are given by 
\begin{align*}
	r_i  = \left(\begin{matrix}0, \cdots,  0, 1,
		r_i^{i+1} ,  \cdots , r_i^{k} ,\cdots, r_i^{N}  \end{matrix}\right) ^T
\end{align*}
with the recursive definition (solving the corresponding lower triangular  linear system)
\begin{align*}
	r_i ^{i+1} =\frac{A_{i+1,i}}{\lambda_i-\lambda_{i+1}}
	\qquad\text{and}\qquad 
	r_i ^{k} =\frac{A_{k,i} + \sum_{j=i+1}^{k-1} A_{k,j} r_i^j}{\lambda_i-\lambda_{k}}
	\text{ for }k=i+2, \ldots, N
	.
\end{align*}
We prove now a simple explicit expression for the eigenvectors, i.e. for $k=i+1, \ldots, N$ we prove
\begin{align*}
	r_i ^{k} =-\frac{f_k}{1- \sum_{j=0}^i f_j}.
\end{align*}
The proof is done by induction over $k$ starting with $k=i+1$.
We have 
\begin{align*}
	r_i ^{i+1} &= \frac{A_{i+1,i}}{\lambda_i -\lambda_{i+1}}\\
	&=  \frac{\alpha (v_{i+1}- v_i) f_{i+1}}{(v_{i+1}- v_i)- \alpha \left(\sum_{j=i+1}^N (v_{j}- v_i)f_j - \sum_{j=i+2}^N (v_{j}- v_{i+1})f_j\right)}\\
	&= - \frac{ f_{i+1}}{\frac{1}{\alpha}+ \frac{1}{v_{i+1}-v_i} \left(\sum_{j=i+2}^N (v_{j}- v_i - (v_j - v_{i+1}))f_j 
		+ (v_{i+1} -v_i ) f_{i+1}\right)}\\ 
	&= - \frac{ f_{i+1}}{1- \rho+  \left(\sum_{j=i+2}^N f_j 
		+  f_{i+1}\right)}\
	=-\frac{f_{i+1}}{1- \sum_{j=0}^i f_j}.
\end{align*}
Assuming now for $k>i+1$
\begin{align*}
	r_i ^{j} =-\frac{f_{j}}{1- \sum_{l=0}^i f_l}, j=i+1, \ldots, k-1,
\end{align*}
we obtain for $k>i+1$
\begin{align*}
	r_i ^{k} &= \frac{A_{k,i}+ \sum_{j=i+1}^{k-1} A_{k,j} r_i^j}{\lambda_i -\lambda_k}\\
	&=  \frac{\alpha (v_{k}- v_i) f_{k} + \alpha \sum_{j=i+1}^{k-1} (v_k-v_j) f_k r_i^j}{(v_{i}- v_k)- \alpha \left(\sum_{j=i+1}^N (v_{j}- v_i)f_j - \sum_{j=k+1}^N (v_{j}- v_{k})f_j\right)}\\
	&=  \frac{ f_{k} \left[  (v_{k}- v_i)  + \sum_{j=i+1}^{k-1} (v_k-v_j)  r_i^j \right]}{\frac{1}{\alpha} (v_i-v_k) -  \left(\sum_{j=i+1}^N (v_{j}- v_i - (v_j - v_{k}))f_j 
		+ \sum_{j=i+1}^k(v_{j} -v_k ) f_{j}\right)}\,
\end{align*}
which gives 
\begin{align*}
	r_i ^{k} &= -  \frac{ f_{k} \left[  (v_{k}- v_i)  + \sum_{j=i+1}^{k-1} (v_k-v_j)  r_i^j \right]}{\left[(1-\rho) (v_k-v_i) -  \sum_{j=i+1}^N (v_{k}- v_i )f_j \right]
		- \sum_{j=i+1}^k(v_{j} -v_k ) f_{j}}\\ 
	&= -  \frac{ f_{k} \left[  (v_{k}- v_i)  + \sum_{j=i+1}^{k-1} (v_k-v_j)  r_i^j \right]}{(1-\sum_{j=0}^i f_j) \left( (v_k-v_i) +  \sum_{j=i+1}^{k-1} (v_{k}- v_j) r_i^j\right) }\\ 
	&=-\frac{f_{k}}{1- \sum_{j=0}^i f_j}.
\end{align*}
This proves the explicit expression.

\subsection{Total linear degeneracy}
Using this representation one can easily prove that the hyperbolic system is totally linear degenerate \cite{CHN09} computing 
$
\nabla \lambda_i \cdot r_i
$
and proving that this is equal to zero for all $i=0, \ldots, N$. We compute
\begin{align*}
	\nabla \lambda_i \cdot r_i = \left(\begin{matrix}\lambda_i-v_i\\ \cdots \\ \lambda_i-v_i \\ \lambda_i-v_{i+1}\\
		\cdots \\ \lambda_i-v_{N}
	\end{matrix}\right)^T \cdot \left(\begin{matrix}0\\ \vdots \\ 0 \\ 1 \\  r_i^{i+1}\\
		\cdots \\ r_i^{N}
	\end{matrix}\right) 
\end{align*}
This yields 
\begin{align*}
	\nabla \lambda_i \cdot r_i &= \lambda_i-v_i  +  \sum_{j=i+1}^N (\lambda_i-v_{j})  r_i^j\\
	&= \lambda_i-v_i  -   \sum_{j=i+1}^N (\lambda_i-v_{j})  \frac{f_j}{1-\sum_{l=0}^i f_l}\\
	&= \frac{1}{1-\sum_{l=0}^i f_l} \left((\lambda_i-v_i ) (1-\sum_{j=0}^i f_j)-   \sum_{j=i+1}^N (\lambda_i-v_{j}) f_j\right).
\end{align*}
This expression  is equal to zero, since
\begin{align*}
	&(\lambda_i-v_i ) (1-\sum_{j=0}^i f_j)-   \sum_{j=i+1}^N (\lambda_i-v_{j}) f_j\\
	&=\lambda_i (1-\sum_{j=0}^N f_j)-   v_i (1 -\sum_{j=0}^i f_j) + \sum_{j+1}^N v_{j} f_j\\
	&=\lambda_i (1-\sum_{j=0}^N f_j)-   v_i (1 -\sum_{j=0}^N f_j) + \sum_{j+1}^N v_i f_j  - \sum_{j+1}^N v_i f_j + \sum_{j+1}^N (v_{j}-v_i) f_j\\
	&=\lambda_i (1-\sum_{j=0}^N f_j)-   v_i (1 -\sum_{j=0}^N f_j) + \sum_{j+1}^N (v_{j}-v_i) f_j
\end{align*}
and this is equal to
\begin{align*}
	&(\lambda_i -v_i)(1-\sum_{j=0}^N f_j)+ \sum_{j+1}^N (v_{j}-v_i) f_j\\
	&=\frac{1}{\alpha} \left((\lambda_i -v_i)+ \alpha \sum_{j+1}^N (v_{j}-v_i) f_j\right)=0
\end{align*}
due to the definition of the eigenvalues $\lambda_i$.
Thus,  we have obtained a totally linear degenerate hyperbolic system.

In the following subsections we determine Riemann invariants for the hyperbolic system and investigate  the geometry of the integral curves.

\subsection{Riemann invariants}

In the present case the Lax curves are given  by the integral curves due to the linear degeneracy of all fields.
To determine them we look first for Riemann invariants  $w_i^j$ such that $\nabla w_i^j  \cdot r_i =0 $ for $i=0, \cdots N$ and $j= 0, \ldots, N, j \neq i$, i.e.
$N$ Riemann invariants for each field.
We claim that the following functions are Riemann invariants
for  each field $i=0, \ldots, N$.
\begin{align*}
	w_i^j = f_j, j=0, \cdots , i-1
\end{align*}
for $i= 1, \ldots, N$.
Moreover, we have the  following  invariants  for $i=0, \ldots, N-1$ given by
\begin{align*}
	w_i^j=  r_i^j , j=i+1, \dots,  N.
\end{align*}

This is obvious for $f_j, j=0, \ldots , i-1$ due to the form of the eigenvectors $r_i$.
The invariance of 
$w_i^j=  r_i^j , j=i+1, \ldots,  N, i=0, \ldots, N-1$ is seen as follows.
Compute
\begin{align*}
	&\nabla w_i^j \cdot r_i = \nabla r_i^j \cdot r_i= - \nabla \left( \frac{f_j}{1- \sum_{l=0}^i f_l} \right) \cdot r_i=\\
	& -\left(\begin{matrix}\overbrace{\frac{f_j}{(1- \sum_{l=0}^i f_l)^2} , \cdots , \frac{f_j}{(1- \sum_{l=0}^i f_l)^2}}^{0, \cdots, i}  ,  0 , \cdots , 0 ,
		\overbrace{\frac{1}{1- \sum_{l=0}^i f_l}}^{j},0
		, \cdots , 0
	\end{matrix}\right) \cdot \left(\begin{matrix}0\\ \vdots \\ 0 \\ 1 \\  r_i^{i+1}\\
		\cdots \\ r_i^{N}
	\end{matrix}\right) .
\end{align*}
This expression is equal to 
\begin{align*}
	& - \frac{f_j}{(1- \sum_{l=0}^i f_l)^2} - \frac{1}{1- \sum_{l=0}^i f_l} r_i^j\\
	&= - \frac{f_j}{(1- \sum_{l=0}^i f_l)^2} + \frac{1}{1- \sum_{l=0}^i f_l} \frac{f_j}{1- \sum_{l=0}^i f_l} =0.
\end{align*}

\setcounter{example}{2}
\begin{example}[revisited]
	We consider the 3-equation case $N=2$:
	The Riemann invariants are for the field 
	$i=2$
	given by 
	$f_0$ and $f_1$.
	For 
	$i=1$ they are given by 
	$f_0$ and $r_1^2$
	and for the field 
	$i=0$ we have 
	$r_0^2$
	and  $ r_0^1 $.
\end{example}

\section{Richness and conservative form  of the discrete velocity model}
\label{discreteinvest2}

We  prove that the system can be written in diagonal and in conservative  form using suitable transformations between
N-simplexes and the unit N-cube.

\subsection{The diagonal form of the  system}

Out of the Riemann invariants computed above, we choose a special subset of $N+1$ Riemann invariants $w_k, k=0, \ldots, N$ such that 
for all $i= 0, \ldots, N $ with $i \neq k$ we have 
\begin{align*}
	\nabla w_k \cdot r_i =0.
\end{align*}

These invariants are obtained by
choosing 
\begin{align*}
	w_0 =  f_0 \in [0,1], \;
	w_k = r_{k-1}^{k} = \frac{f_{k}}{1- \sum_{j=0}^{k-1} f_j} \in [0,1], k= 1, \ldots, N
\end{align*}
Note that the transformation  from $f_i$ to $w_i$ variables is a transformation from the  unit N-simplex to the unit N-cube.

We prove that for all $i= 0, \ldots, N $ with $i \neq k$ we have 
$
\nabla w_k \cdot r_i =0.
$
A standard computation gives then that the transformed system is in diagonal form.

For $w_0$ the statement is obvious.

For $k=1, \ldots, N$ and $i=0, \ldots, N, i \neq k$
\begin{align*}
	\nabla w_k \cdot r_i& = \nabla r_{k-1}^k \cdot r_i= - \nabla \left( \frac{f_k}{1- \sum_{l=0}^{k-1} f_l} \right) \cdot r_i\\
	&= -\left(\begin{matrix}\overbrace{\frac{f_k}{(1- \sum_{l=0}^{k-1} f_l)^2} , \cdots , \frac{f_k}{(1- \sum_{l=0}^{k-1} f_l)^2} }^{0, \cdots, k-1},  
		\frac{1}{1- \sum_{l=0}^{k-1} f_l} ,0 , \cdots , 0 
	\end{matrix}\right) \cdot \left(\begin{matrix}0\\ \vdots \\ 0 \\ 1 \\  r_i^{i+1}\\
		\cdots \\ r_i^{N}
	\end{matrix}\right) .
\end{align*}
This is obviously equal to zero, as long as $i>k$. In case  $i<k$ it is equal to
\begin{align*}
	&  \frac{f_k}{(1- \sum_{l=0}^{k-1} f_l)^2} \left( 1+ r_i^{i+1} + \cdots +r_i^{k-1}\right) + \frac{1}{1-\sum _{j=0}^{k-1} f_j} r_i^k.
	\\
	&= \frac{1}{(1- \sum_{l=0}^{k-1} f_l)^2} \left(f_k(1+ \sum_{j=i+1}^{k-1} r_i^{j} ) +(1-\sum _{j=0}^{k-1} f_j) r_i^{k} \right)
	\\
	&= \frac{1}{(1- \sum_{l=0}^{k-1} f_l)^2 (1-\sum_{l=0}^{i} f_l)} \left(f_k(1-\sum_{j=0}^{i} f_j- \sum_{j=i+1}^{k-1} f_{j}  ) -(1-\sum _{j=0}^{k-1} f_j) f_k\right)=0.
\end{align*}
This leads to the equations in diagonal form
\begin{align}
	\partial_t  w_i + \lambda_i  (w)  \partial_x w_i= 0 ,  \; i=0, \ldots,  N .
\end{align}
with
\begin{align*}
	\lambda_0= - \frac{1}{1-\rho} \sum_{j=1}^{N} v_j w_j \prod_{l=0}^{j-1}(1-w_l)   \le  0, \; \lambda_N = 1 \\
	\lambda_i = v_i- \frac{1}{1-\rho} \sum_{j=i+1}^{N} (v_j-v_i) w_j \prod_{l=0}^{j-1}(1-w_l)  , i=1, \ldots, N-1
\end{align*}

\begin{remark}
	Note that the integral curves in the  transformed variables are simply given by straight lines in   the direction of the coordinate axes. Therefore the unit cube is obviously an invariant domain for the transformed equations.
	In primitive  variables $f$  this gives  the invariance of the unit N-simplex for the equations.
\end{remark}

\begin{remark}
	Note that the backwards transformation is given by
	\begin{align*}
		f_0 =  w_0, \; 
		f_k = w_k \prod_{j=0}^{k-1}(1-w_j),\;
		1-\rho = \prod_{j=0}^{N}(1-w_j).
	\end{align*}
\end{remark}

\setcounter{example}{2}
\begin{example}[revisited]
	We consider again  the 3-equation case $N=2$. With the 
	invariants $f_0$ and $r_0^1$ for the field $i=2$, the 
	invariants $f_0$ and $r_1^2$ for $i=1$ and 
	invariants $r_0^1$ and $r_1^2$ for $i=0$.
	The eigenvalues in the new variables are 
	\begin{align*}
		\lambda_0  (w)  =  - \frac{v_1w_1+w_2(1-w_1)}{(1-w_1)(1-w_2))}\\
		\lambda_1  (w) = v_1 - (1-v_1) \frac{w_2}{1-w_2} =  \frac{v_1 - w_2}{1-w_2} \\
		\; \lambda_2  (w) =1
	\end{align*}
	and the macroscopic quantities are 
	\begin{align*}
		\rho = 1 -  (1-w_0)(1-w_1)(1-w_2)\\
		q=v_1 f_1 +f_2 = (1-w_0) (v_1 w_1  + w_2 - w_1 w_2).
	\end{align*}
\end{example}

In a last step  we determine a second transformation such that the equations can be written in conservative from.
Hyperbolic systems which can be written in diagonal and conservative form are called rich \cite{serre,peng}.

\subsection{Rich systems and the conservative form}

We define the functions $N_k (w) , k=0, \ldots, N$ as functions of the above Riemann invariants $w_k$ as 
\begin{align*}
	N_k (w) = \prod_{j=k}^{N}(1-w_j) \in [0,1].
\end{align*}
The transformation from $w_i$ to $N_i$ variables  is a transformation from the unit N-cube into  an N-simplex.
Note that for $k=0, \ldots, N-1$
\begin{align*}
	\lambda_k (w) &= v_k- \frac{1}{1-\rho} \sum_{j=k+1}^{N} (v_j-v_k) f_j \\
	&=v_k- \frac{ \sum_{j=k+1}^{N} (v_j-v_k) w_j \prod_{l=0}^{j-1}(1-w_l)}{\prod_{j=0}^{N}(1-w_j)}\\
	&=v_k - \frac{ \sum_{j=k+1}^{N} (v_j-v_k) w_j \prod_{l=0}^{j-1}(1-w_l)}{\prod_{j=0}^{N}(1-w_j)}\\
	&=v_k - \frac{  \sum_{j=k+1}^{N} (v_j-v_k) w_j \prod_{l=k}^{j-1}(1-w_l)}{\prod_{j=k}^{N}(1-w_j)}\\
	&=v_k - \frac{ (v_{k+1}-v_k) w_{k+1}+ \sum_{j=k+2}^{N} (v_j-v_k) w_j \prod_{l=k+1}^{j-1}(1-w_l)}{\prod_{j=k+1}^{N}(1-w_j)}
\end{align*}
and therefore
\begin{align*}
	\frac{\partial \lambda_k (w) }{\partial w_j } =0
\end{align*}
for $j \le  k$.
For $j > k$ we have 
\begin{align}
	\label{lamex}
	\frac{\partial \lambda_k (w) }{\partial w_j } = - \frac{\lambda_j - \lambda_k}{1-w_j}.
\end{align}
A proof can be found in the appendix.

Moreover, we have directly 
\begin{align*}
	\frac{\partial N_k (w) }{\partial w_j }=\frac{\partial  }{\partial w_j } (\prod_{l=k}^{N}(1-w_l))
	=- \prod_{l=k , l \neq  j}^{N}(1-w_l)
\end{align*}
if $j\ge k$ and $ \frac{\partial N_k (w) }{\partial w_j }=0$ if  $j < k$.
Altogether,  we  have
\begin{align}
	\label{pengcrit}
	(\lambda_j - \lambda_k) \frac{\partial N_k (w) }{\partial w_j }= N_k(w) \frac{\partial \lambda_k (w) }{\partial w_j }.
\end{align}
This is obvious for $j \le  k$. For $j >  k$ we have 
\begin{align*}
	(\lambda_j - \lambda_k) \frac{\partial N_k (w) }{\partial w_j }
	&=-  (\lambda_j - \lambda_k)  \prod_{l=k , l \neq  j}^{N}(1-w_l)\\
	&= (1-w_j)\prod_{l=k,l \neq  j}^{N}(1-w_l) \frac{\partial \lambda_k (w) }{\partial w_j }
	= \prod_{l=k}^{N}(1-w_l) \frac{\partial \lambda_k (w) }{\partial w_j }\\
	&= N_k(w) \frac{\partial \lambda_k (w) }{\partial w_j }.
\end{align*}
Using (\ref{pengcrit}) a classical computation, see \cite{peng,serre}, gives  the conservative form
\begin{align*}
	\partial_t N_i(w) + \partial_x (\lambda_i(w) N_i(w))=0, i=0, \ldots, N.
\end{align*}
The eigenvalues in the new variables are 
\begin{align*}
	\lambda_0= - \frac{1}{1-\rho} \left(\sum_{j=1}^{N-1} v_j (1-\frac{N_j}{N_{j+1}}) \prod_{l=0}^{j-1} \frac{N_l}{N_{l+1}}+  (1-N_N) \prod_{l=0}^{N-1}\frac{N_l}{N_{l+1}}\right)  
\end{align*}
and
for $i=1, \ldots, N-2$
\begin{align*}
	\lambda_i = v_i- \frac{1}{1-\rho} \left(\sum_{j=i+1}^{N-1} (v_j-v_i) (1-\frac{N_j}{N_{j+1}}) \prod_{l=0}^{j-1} \frac{N_l}{N_{l+1}} + (1-v_i) (1-N_N) \prod_{l=0}^{N-1}\frac{N_l}{N_{l+1}}   \right)
\end{align*}
and
\begin{align*}
	\lambda_{N-1} = v_{N-1}- \frac{1}{1-\rho}  (1-v_{N-1}) N_0 (\frac{1}{N_N}-1) 
	,\qquad \lambda_{N}=1
	.
\end{align*}
With $N_0 = 1-\rho$ this can be rewritten as 
\begin{align*}
	\lambda_0= -  \left(\sum_{j=1}^{N-1} v_j (\frac{1}{N_j}-\frac{1}{N_{j+1}})  +  (\frac{1}{N_N}-1) \right)  
\end{align*}
and
for $i=1, \ldots, N-2$
\begin{align*}
	\lambda_i = v_i- \left(\sum_{j=i+1}^{N-1} (v_j-v_i) (\frac{1}{N_{j}}- \frac{1}{N_{j+1}} )+ (1-v_i) ( \frac{1}{N_{N}}-1)   \right)
\end{align*}
and
\begin{align*}
	\lambda_{N-1} = v_{N-1}-  (1-v_{N-1}) (\frac{1}{N_N}-1) 
	,\qquad \lambda_{N}=1
	.
\end{align*}

\begin{remark}
	The backwards transformation is 
	\begin{align*}
		w_N = 1-N_N,\qquad
		w_i = 1- \frac{N_i}{N_{i+1}}, i= 0, \ldots , N-1 
	\end{align*}
	and the direct transformation between primitive and conservative variables is given by  
	\begin{align*}
		f_N = N_0 (\frac{1}{N_N}-1),\qquad
		f_k = N_0( \frac{1}{N_{k}} - \frac{1}{N_{k+1}}), k= 0, \ldots , N-1 \nonumber
	\end{align*}
	and 
	\begin{align*}
		N_0= 1-\rho,\qquad
		N_k = \frac{1- \rho}{1-\sum_{j=0}^{k-1}f_j}, k=1, \ldots , N.
	\end{align*}
\end{remark}

\setcounter{example}{1}
\begin{example}[revisited]
	For the 2-equation case  $N=1$ we have 
	\begin{align*}
		N_0= 1-\rho, \qquad  
		N_1 = \frac{1- \rho}{1-f_0}.
	\end{align*}
	and
	\begin{align*}
		f_1 = N_0 (\frac{1}{N_1}-1),\qquad
		f_0 = N_0( \frac{1}{N_{0}} - \frac{1}{N_{1}}).
	\end{align*}
	The equations are with $	\lambda_0 = 1-  \frac{1}{N_{1}}$ given by 
	\begin{align*}
		\partial_t N_0 + \partial_x (\lambda_0(N) N_0)&=0,\\
		\partial_t N_1 + \partial_x N_1&=0.
	\end{align*}
	Moreover, the flux is given by  $q=f_1 = N_0 (\frac{1}{N_1}-1) $.
\end{example}

\begin{example}[revisited]
	For  the 3-equation case $N=2$ we have 
	\begin{align*}
		N_0 = 1-\rho,\qquad
		N_1 =  (1-w_1) (1-w_2), \qquad
		N_2 =  (1-w_2) 
	\end{align*}
	and 
	\begin{align*}
		w_0 = 1-\frac{N_0}{N_1}, \qquad 
		w_1 =  1- \frac{N_1}{N_2},\qquad
		w_2 =  1-N_2.
	\end{align*}
	This gives with 
	\begin{align*}
		\lambda_0  (N)  &= - \frac{1}{N_1N_2} (v_1(N_2-N_1) +N_1(1-N_2)), \\
		\lambda_1  (N) &= \frac{1}{N_2} (N_2-(1-v_1))
	\end{align*}
	the system of equations
	\begin{align*}
		\partial_t  N_0 +  \partial_x ( \lambda_0(N)N_0)&= 0 \\
		\partial_t  N_1 +   \partial_x ( \lambda_1(N)N_1)&= 0\\
		\partial_t  N_2 +   \partial_x N_2 &= 0.
	\end{align*}
	Moreover, the flux is $q= N_0((1-v_1)\frac{1}{N_2} +  \frac{v_1}{N_1}-1)$.
\end{example}

\begin{remark}
	Defining a continuous version of the transformation between  primitive and conservative variables 
	we have 
	\begin{align*}
		N= \frac{1-\rho}{1-\int_0^v f(v^\prime) d v^\prime}
	\end{align*}
	and
	\begin{align*}
		f= (1-\rho)\frac{N^\prime (v)}{N^2(v)}.
	\end{align*}
	Then,  the equation
	\begin{align}
		\partial_t f + v \partial_x f 
		- \frac{ 1 }{1-\rho} \left[
		\int_v^1 ( \hat v -v )  f(\hat v)  d \hat v  \partial_x f (v)
		- \int_0^v  ( v-\hat v  )   \partial_x f (\hat v)  d \hat v\  f( v)  \right]=0
	\end{align}
	is transformed into the  problem
	\begin{align}
		\label{conseq}
		\partial_t N + \partial_x \left(N\left(1-\int_v^1 \frac{1}{N(v^\prime)} d v^\prime \right) \right) =0 .
	\end{align}
	
\end{remark}

\begin{remark}
	The fact, that the hyperbolic system is rich and totally linear degenerate yields different further properties of the system.
	For example, for the homogeneous system one obtains explicit representation formulas and existence of global entropy solutions,
	see \cite{RS,Li}. For the relaxation system  the so-called semi-linear behaviour is obtained in \cite{CHN09} preventing the 
	appearance of shock solutions.
\end{remark}

\section{The full relaxation model}
\label{relax}

In this section we consider the full relaxation model. In particular, a stability condition is formulated based on the Chapman-Enskog procedure, see \cite{G04,GKR03} for similiar investigations for the ARZ equations.
The relation of this condition to the appearance of persistent periodic solutions will be numerically  investigated in the following Section \ref{numericalresults}.

\subsection{Stability of the continuous kinetic model}

Considering the continuous inhomogeneous problem (\ref{kinbasic}), multiplying the equation with $1$ and $v$  and integrating with respect to $v$ one obtains the balance equations, that is equations for the zeroth and first  moment of the distribution function $f$, i.e.
\begin{align*}
	&\partial_t \rho + \partial_x q =0\\
	&\partial_t q +\partial_x e  +
	\frac{1}{1-\rho } \int \int_{\hat v < v}  ( v-\hat v  )^2 f(v)    \partial_x f (\hat v)  d \hat v\  d v  = - \frac{1}{\epsilon } (q-F(\rho))
\end{align*} 
using the notation $e= \int v^2 f(v) dv $.
Using the fact that  up to order  $\mathcal{O} (\epsilon)$ we have $f = f^e(\rho)$  we obtain up to  order $\mathcal{O} (\epsilon^2)$ 
\begin{align*}
	q = F(\rho)-  \epsilon \Big(E^\prime(\rho) -(F^\prime(\rho))^2+
	\frac{1}{1-\rho } \int \int_{\hat v < v}  ( v-\hat v  )^2 f^e(\rho,v)    \partial_\rho f^e (\rho,\hat v)  d \hat v\  d v \Big) \partial_x \rho
\end{align*} 
with $E=E(\rho) =  \int v^2 f^e(\rho,v) dv$. 
Inserting this into the continuity equation gives a drift-diffusion equation for the density with positive diffusion coefficient 
$D(\rho)$ as long as the  formal Chapman-Enskog stability condition 
\begin{align}
	\label{contstability}
	D(\rho) = -( F^\prime (\rho))^2 + E^\prime (\rho) +
	\frac{1}{1-\rho } \int \int_{\hat v < v}  ( v-\hat v  )^2 f^e( \rho,v)    \partial_\rho f^e (\rho, \hat v)  d \hat v\  d v  \ge 0
\end{align} 
is fulfilled.

\setcounter{example}{0}
\begin{example}[revisited]
	We reconsider Example \ref{excont} with the equilibrium functions $f^e$ depending only on the first 3 moments $\rho,F,E$ and given by 
	\begin{align*}
		f^e(\rho,v) = \delta_0 (v)  \left(\rho-F - (\gamma-1) \frac{F-E}{1-\lambda} \right) +  \frac{F-E}{1-\lambda} \frac{g(v)}{v}+ \delta_1(v) \left(  F-\frac{F-E}{1-\lambda}\right).
	\end{align*}
	We had $g \ge 0$, $\int_0^1  g(v) dv =1$. Moreover,  $\lambda = \int_0^1 v g(v) dv$ and 
	$\gamma = \int_0^1 \frac{g(v)}{v} d v < \infty$.
	Positivity of $f^e$ did require $F- \frac{1-\lambda}{\gamma-1}  (\rho-F) \le E $   and $\lambda F \le E \le F$
	for all $0 < \rho < 1$.

	Then
	\begin{align*}
		&\int \int_{\hat v < v}  ( v-\hat v  )^2 f^e( \rho,v)    \partial_\rho f^e (\rho, \hat v)  d \hat v\  d v \\
		&=  \int_0^1 v  g(v)dv \; \frac{F-E}{1-\lambda}  \left(1-F^\prime - (F^\prime-E^\prime)\frac{\gamma-1}{1-\lambda} \right) \\
		&+\frac{1}{(1-\lambda)^2}  \int \int_{\hat v < v}  \frac{( v-\hat v  )^2}{v \hat v}  g(v) g(\hat v) d \hat v dv   \; (F-E)  (F^\prime-E^\prime)\\
		&
		+   \left(  F-\frac{F-E}{1-\lambda}\right)  \left(1-F^\prime - \frac{\gamma-1}{1-\lambda} (F^\prime-E^\prime) \right)\\
		&+   \int_0^1 \frac{( 1-\hat v  )^2}{\hat v }  g(\hat v) d \hat v  \frac{F^\prime-E^\prime}{1-\lambda} \left(  F-\frac{F-E}{1-\lambda}\right).
	\end{align*}
	Note that 
	$$\int \int_{\hat v < v}  (\frac{v}{ \hat v} + \frac{\hat v }{  v} -2)g(v) g(\hat v) d \hat v dv = \gamma \lambda-1.$$
	Therefore the above is equal to 
	\begin{align*}
		& \frac{\lambda}{1-\lambda} (F-E)  \left(1-F^\prime -  (F^\prime-E^\prime)\frac{\gamma-1}{1-\lambda} \right) \\
		&+ (\frac{\gamma \lambda-1}{(1-\lambda)^2})  (F-E)  (F^\prime-E^\prime)\\
		&
		+   \left(  F-\frac{F-E}{1-\lambda}\right)  \left(1-F^\prime -  (F^\prime-E^\prime) \frac{\gamma-1}{1-\lambda}\right)\\
		&+ \frac{\gamma-2+\lambda}{1-\lambda} (F^\prime-E^\prime)\left(  F-\frac{F-E}{1-\lambda}\right)\\
		&= E - (F^\prime -E^\prime)F -F^\prime E 
	\end{align*}
	Thus, we obtain the stability condition
	\begin{align}
		\label{condstab}
		-(F^\prime)^2 + E^\prime  + \frac{1}{1-\rho} (E - (F^\prime -E^\prime)F -F^\prime E ) \ge 0.
	\end{align}
	
	In case 
	$E=F$ or $
	f^e(\rho,v) = \delta_0 (v)  \left(\rho-F  \right) +\delta_1(v)   F$ 
	this is simply 
	\begin{align*}
		(1- F^\prime )  ( F^\prime + \frac{F}{1-\rho} ) \ge 0,
	\end{align*}
	which is fulfilled for all concave $F$.
	In case $E \ne F$ unstable regions for $\rho$ can be easily obtained even for concave first moments $F$, see the Section \ref{numericalresults}.
\end{example}

\subsection{Stability of the discrete velocity models with relaxation}

For the full discrete velocity system with relaxation (\ref{discretekin}) with $N+1$ equations, the Chapman Enskog stability condition translates to
\begin{align*}
	-(\partial_\rho F(\rho))^2 + \partial_\rho E (\rho) +
	\frac{1}{1-\rho } \sum_{i=1}^N \sum_{j= 0}^{i-1}( v_i- v_j  )^2 f^e_i( \rho)    \partial_\rho f^e_j (\rho)  \ge 0,
\end{align*} 
where $$ 
F(\rho) = \sum_{i=0}^N v_i f^e_i (\rho) , \;  E(\rho) = \sum_{j=1}^N v_j^2 f^e_i(\rho). $$

\begin{remark}
	Since $\lambda_0 = -\frac{q}{1-\rho}$ 
	the sub characteristic condition \cite{Liu87} for the system  is  given by 
	$$
	-  \frac{ F(\rho)}{1-\rho} = \lambda_0 \le F^\prime(\rho) \le 1 = \lambda_N \  \mbox{ for } \  0 \le \rho \le 1 \ .
	$$
	This  is fulfilled for example for concave fundamental diagrams $F$.
\end{remark}

\setcounter{example}{1}
\begin{example}[revisited]
	$N=1$.
	Sub characteristic and Chapman Enskog stability condition are equivalent and are given by 
	\begin{align*} 
		\begin{aligned}
			(F^\prime (\rho) + \frac{F(\rho)}{1-\rho} ) \left(1 -  F^\prime  (\rho)  \right)
			\ge 0.
		\end{aligned}
	\end{align*}

\end{example}

\begin{example}[revisited]
	$N=2$.
	In this case the  Chapman Enskog stability condition is as in the above continuous example
	given by 
	\begin{align*}
		\begin{aligned}
			-(F^\prime)^2 + E^\prime  + \frac{1}{1-\rho} \left[E -(E+F)   F^\prime   +F E^\prime  \right]
			\ge 0.
		\end{aligned}
	\end{align*}
	Note that $E=F$ gives again the sub-characteristic condition.
	Moreover, realizability requires  functions $F$ and $E$   to lie in the realizability region
	\begin{align*}
		F- v_1(\rho-F) \le E,  \; v_1 F\le E \le F.
	\end{align*}
\end{example}

\begin{example}[revisited]
	
	Realizability requires again $F- \frac{1-\lambda}{\gamma-1}  (\rho-F) \le E$  and $\lambda F \le E \le F$ for all $0 < \rho < 1$ as in the continuous case.
	The stability condition is the same as in \eqref{condstab}.
	
\end{example}
\section{Numerical Method and Results}
\label{numericalresults}

The results obtained in the previous sections can be used to obtain a 
mass-conserving and asymptotic preserving algorithm respecting the invariant domain $\Delta_N$ of the equations.

\subsection{The Godunov method}

Using the transformation to diagonal form the Godunov method can be realized in a particularly efficient way.
Considering the conservative form 
\begin{align*}
	\partial_t N_i(w) + \partial_x (\lambda_i(w) N_i(w))=0, i=0, \ldots, N
\end{align*}
we have  to determine the solution of the Riemann problem for any combination of left and  right states $N^L$ and  $N^R$.

Instead of doing this directly in the conservative variables we change to the diagonal variables
\begin{align*}
	w_N = 1-N_N,\qquad
	w_i = 1- \frac{N_i}{N_{i+1}},\  i= 0, \ldots , N-1.
\end{align*}
The $N$ intermediate states of the solution of the Riemann problem are then simply
\begin{align*}
	w^{M,l}= (w_0^R , \ldots , w_l^R,w^L_{l+1}, \ldots, w^L_N),  l= 0,  \;\ldots , N-1.
\end{align*}

Next, we determine the eigenvalues $\lambda_l, l= 1, \ldots,N-1$ for the intermediate states.
Finally, we determine   the  intermediate state  $w^M = w^{M,l_0}$ such that  
$\lambda_{l_0} (w^{M})<0$ and $\lambda_{l_0+1}(w^{M})>0$. This is the Godunov state for our Riemann problem
Note that we have always $\lambda_0 <0$ and $\lambda_N =1>0$. 
Finally, the transformation 
\begin{align*}
	N^M_k (w) = \prod_{j=k}^{N}(1-w_j^{M}) \in [0,1].
\end{align*}
gives the new Godunov  state in conservative variables  for the Godunov method.

\subsection{Full numerical algorithm}

We use a splitting scheme solving first the advection and then the relaxation part of the equation.
Suppose we know $f_i, i=0, \ldots , N$ at the initial time. 

\noindent
{\em Step 1:} Determine the conservative variables
\begin{align*}
	N_0= 1-\rho,\qquad
	N_k = \frac{1- \rho}{1-\sum_{j=0}^{k-1}f_j},\quad  k=1, \ldots , N.
\end{align*}
{\em Step 2: } Use an e-consistent first order scheme, for example, the above Godunov method,
to solve  for one time-step the conservative advection step 
\begin{align*}
	\partial_t N_i + \partial_x (\lambda_i(N) N_i)=0,\quad  i=0, \ldots, N
\end{align*}
with the eigenvalues $\lambda_N =1$
and
for $i=0, \ldots, N-1$
\begin{align*}
	\lambda_i &= v_i- \sum_{j=i+1}^{N} (v_j-v_i) (\frac{1}{N_{j}}- \frac{1}{N_{j+1}} ) 
	= v_i- \sum_{j=i+1}^{N} (v_j-v_{j-1}) \frac{1}{N_{j}}+(1-v_i) \\
	&= 1- \sum_{j=i+1}^{N} (v_j-v_{j-1}) \frac{1}{N_{j}},
\end{align*}
where we have defined $N_{N+1}=1$.

\noindent
{\em Step 3:} 
Determine the new variables 
\begin{align*}
	f_N = N_0 (\frac{1}{N_N}-1),\qquad 
	f_k = N_0( \frac{1}{N_{k}} - \frac{1}{N_{k+1}}),\quad k= 0, \ldots , N-1.
\end{align*}
{\em Step 4: }
Solve the relaxation step  with an  implicit Euler scheme 
\begin{align*}
	\partial_t f_i =- \frac{1}{T} (f_i-f^e_i(\rho))\quad i=0, \ldots,N.
\end{align*}

Higher order numerical methods can be developed along the usual lines for the numerical treatment of hyperbolic relaxation systems,
see, for example,  \cite{RP05}.

We consider now several numerical examples to illustrate the behaviour of the equations.
In  all cases the computational domain is $[0,1]$. In test-case 1 free boundary conditions are used. In the other test-cases
we use periodic boundary conditions.
We chose for all examples a fine spatial discretization with $2000$ discretization points and a corresponding time discretization resulting from the CFL condition.

\subsection{Test-case 1: Riemann problems}

As a first example, we consider a Riemann problem for the hyperbolic equations without relaxation term.
We choose $N=10$ and  4 different   kinetic  initial conditions. They are always chosen such that for $x \le \frac{1}{2}$ we have $\rho  = \rho_L= 0.6$ and for $x>\frac{1}{2}$ we have $\rho =  \rho_R=0.8$. In all cases  $q_L$ and $q_R$ are given by $q=F(\rho)$ with $F = \rho(1-\rho) $.

The  first kinetic  initial condition is obtained using the left and right  initial distributions $f_L = \alpha_L \delta_{0.4} + \beta_L \delta_0$ and $f_R = \alpha_R \delta_{0.4} + \beta_R \delta_{0.2}$, where 
$\alpha_{L/R}$ and $\beta_{L/R}$ are chosen such that the above macroscopic initial values are obtained.
The second initial condition is given by using the same $f_L$ and changing the values for  $f_R$ into $f_R=\alpha_R \delta_{0.4} + \beta_R \delta_{0.1}$
The third condition uses  $f_{L/R}=\alpha_{L/R} \delta_{0.2} + \beta_{L/R} \delta_{0}$ and the fourth one
$f_{L/R}=\alpha_{L/R} \delta_{0.4} + \beta_{L/R}  \delta_{0}$ .

\begin{figure}[h!]
	\tikzsetnextfilename{Test1_1}
		\begin{tikzpicture}[scale=0.7]
			\begin{axis}[ylabel = $\rho$,xlabel =  $x$,
				legend style = {at={(0.5,1)},xshift=0.0cm,yshift=0.1cm,anchor=south},
				legend columns= 4,
				]
				\addplot[color = red, thick] file{datariemann/figrho1.txt};
				\addlegendentry{IC 1}				
				\addplot[color = green, thick] file{datariemann/figrho2.txt};
				\addlegendentry{IC 2}
				\addplot[black,thick] coordinates {(0,0.6)  (0.34,0.6)};
				\addplot[black,thick] coordinates {(0.34,0.8)  (1,0.8)};
				\addplot[black,thick] coordinates {(0.34,0.6)  (0.34,0.8)};	
				\addlegendentry{LWR}	
			\end{axis}
		\end{tikzpicture}
	\tikzsetnextfilename{Test1_2}
		\begin{tikzpicture}[scale=0.7]
			\begin{axis}[ylabel = $q$,xlabel =  $x$,
				legend style = {at={(0.5,1)},xshift=0.0cm,yshift=0.1cm,anchor=south},
				legend columns= 4,
				]
				
				\addplot[color = red,thick] file{datariemann/figq1.txt};
				\addlegendentry{IC 1}
				\addplot[color = green,thick] file{datariemann/figq2.txt};
				\addlegendentry{IC 2}
				\addplot[black,thick] coordinates {(0,0.24)  (0.34,0.24)};
				\addplot[black,thick] coordinates {(0.34,0.16)  (1,0.16)};
				\addplot[black,thick] coordinates {(0.34,0.24)  (0.34,0.16)};	
				\addlegendentry{LWR}	
			\end{axis}
		\end{tikzpicture}
	
	\tikzsetnextfilename{Test1_3}
		\begin{tikzpicture}[scale=0.7]
			\begin{axis}[ylabel = $\rho$,xlabel =  $x$,
				legend style = {at={(0.5,1)},xshift=0.0cm,yshift=0.1cm,anchor=south},
				legend columns= 4,
				]
				\addplot[color = red, thick] file{datariemann/figrho3.txt};
				\addlegendentry{IC 3}				
				\addplot[color = blue, thick] file{datariemann/figrho4.txt};
				\addlegendentry{IC 4}
				\addplot[color = black, thick] file{datariemann/figrhorascle.txt};
				\addlegendentry{ARZ}
			\end{axis}
		\end{tikzpicture}
	\tikzsetnextfilename{Test1_4}
		\begin{tikzpicture}[scale=0.7]
			\begin{axis}[ylabel = $q$,xlabel =  $x$,
				legend style = {at={(0.5,1)},xshift=0.0cm,yshift=0.1cm,anchor=south},
				legend columns= 4,
				]
				
				\addplot[color = red,thick] file{datariemann/figq3.txt};
				\addlegendentry{IC 3}
				\addplot[color = blue,thick] file{datariemann/figq4.txt};
				\addlegendentry{IC 4}
				\addplot[color = black,thick] file{datariemann/figqrascle.txt};
				\addlegendentry{ARZ}
			\end{axis}
		\end{tikzpicture}
	
	\caption{Density $\rho$ and flux $q$ for Riemann problems with $\rho_L= 0.6, \rho_R=0.8$, $q =\rho(1-\rho)$ and 4 different kinetic initial situations. For comparison LWR and ARZ solutions are shown respectively.}
	\label{fig1}
\end{figure}
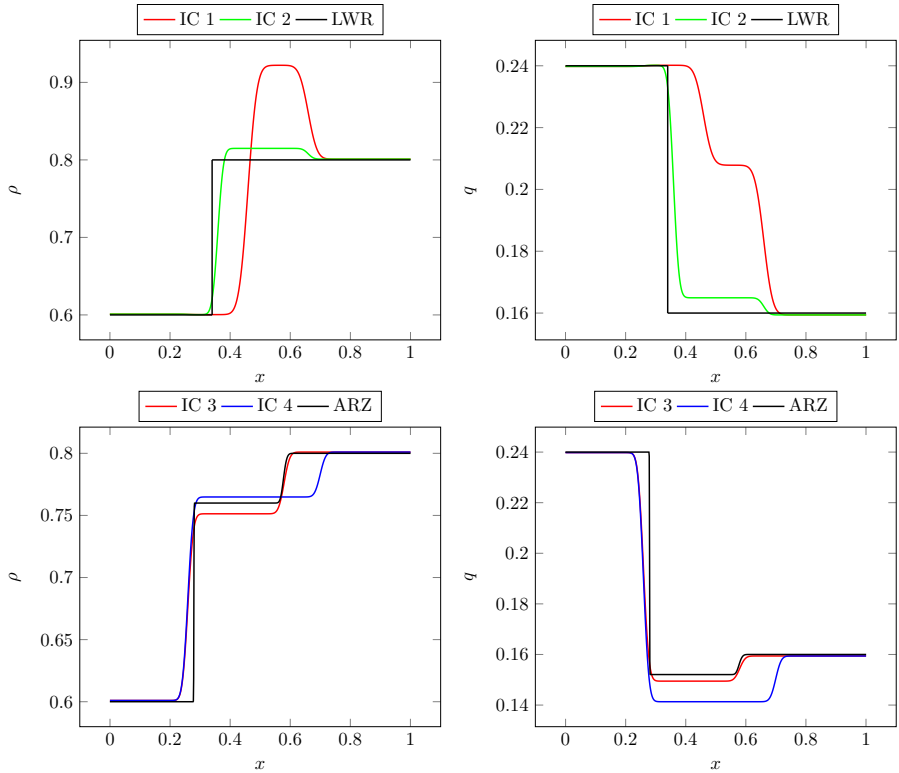

Figure \ref{fig1} top row shows the solutions of the Riemann problems at $t=0.4$ for initial states 1 and 2 and the solution of the LWR equations with
flux function $ \rho(1-\rho) $.  Figure \ref{fig1} bottom  row shows the  solutions of the Riemann problems  for initial states 3 and4 4 and the solution of the ARZ equations $\partial_t \rho+ \partial_x (\rho u) =0$, $\partial_t u + u \partial_x u + \rho p^\prime (\rho) \partial_x u =0 $ without right hand side and pressure $p$ given by
$p(\rho) = \frac{\rho}{1-\rho}$.

One observes in these figures a variety of solutions for the same $\rho$ and $q$ initial values depending on the kinetic initial conditions.

\subsection{Test-case 2:  Decay to equilibrium for stable situation}

We reconsider Example \ref{exN2} and the full relaxation model with $\epsilon= 0.01$. 
As a stable  example,  consider $E= F(\rho)(1-\alpha \rho) $.
Then, the stability condition is 
\begin{align*}
	\begin{aligned}
		-&(F^\prime(\rho))^2 + E^\prime (\rho) + \frac{1}{1-\rho} \left[E (\rho) -(F+E)  (\rho) F^\prime  (\rho) +F  (\rho) E^\prime (\rho) \right]\\
		& =  (1-F^\prime) \left(F^\prime + \frac{F}{1-\rho} \right)
		- \alpha \left(
		F^\prime   \rho+ F + \frac{F}{1-\rho} \left[ \rho+ F ) \right] \right) \ge 0.
	\end{aligned}
\end{align*}
Choosing the concave LWR fundamental diagram $F(\rho ) = \rho(1-\rho)$, the above expression is 
$ 2 \rho (1- \rho)
- \alpha\left(
\rho(2+\rho ) (1-\rho) \right)$.
This is larger than zero if 
\begin{align*}
	\alpha \le  \frac{2}{3} \le  \frac{2}{2+\rho}.
\end{align*}
In case $v_1=\frac{1}{2}$, the function $E= F(\rho)(1-\alpha \rho) $ fulfils  the realizability constraint (\ref{real}) as long as 
$\alpha \le \frac{1}{2}$. 

For the numerical experiments we consider  $E=F$ and  $E=F(1-\frac{\rho}{2})$.
The  initial conditions are defined as
$$
f_i (0,x) =  \frac{1}{3}(0.7+ 0.1\sin(6 \pi x)),\quad i= 0,1,2.
$$

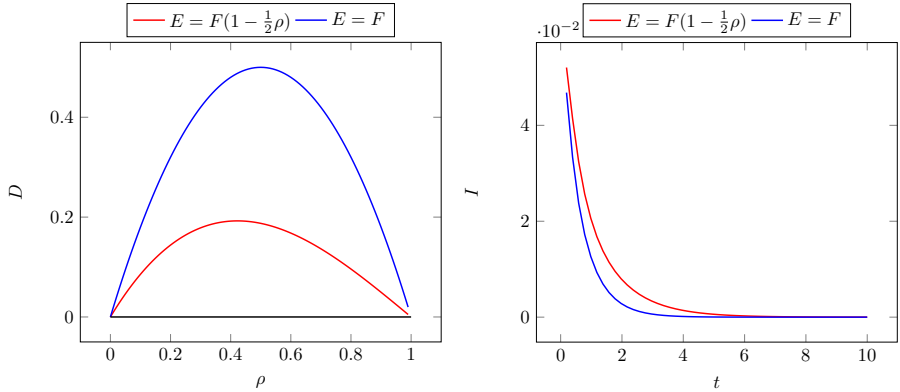
\begin{figure}[h!]
	\tikzsetnextfilename{stability}
		\begin{tikzpicture}[scale=0.7]
			\begin{axis}[ylabel = $D$,xlabel =  $\rho$,
				legend style = {at={(0.5,1)},xshift=0.0cm,yshift=0.1cm,anchor=south},
				legend columns= 4,
				]
				\addplot[color = red,thick] file{datafunctions/stabcond1.txt};
				\addlegendentry{$E= F(1-\frac{1}{2}\rho)$}
				\addplot[color = blue,thick] file{datafunctions/stabcond2.txt};
				\addlegendentry{$E=F$}
				\addplot[black,thick] coordinates {(0,0.0)  (1,0.0)};
			\end{axis}
		\end{tikzpicture}
	\tikzsetnextfilename{decay}
		\begin{tikzpicture}[scale=0.7]
			\begin{axis}[ylabel = $I$,xlabel =  $t$,
				legend style = {at={(0.5,1)},xshift=0.0cm,yshift=0.1cm,anchor=south},
				legend columns= 4,
				]
				\addplot[color = red,thick] file{data/stable1.txt};
				\addlegendentry{$E= F(1-\frac{1}{2}\rho)$}
				\addplot[color = blue,thick] file{data/stable2.txt};
				\addlegendentry{$E=F$}
			\end{axis}
		\end{tikzpicture}
	\caption{Decay to constant state for for $N=2$ with $F= \rho(1-\rho)$  and $E= F(1-\frac{1}{2}\rho)$, $E=F$, respectively. Left:  plot of 
		stability conditions. Right:  plot of decay rate $I$ versus time.  }
	\label{figdecay}
\end{figure}

In Figure \ref{figdecay}   the Chapman-Enskog stability term  $D(\rho)$ is plotted on the left.
Figure \ref{figdecay}  on the right shows a plot of the quantity
$$
I (t)= \int_0^1 \vert \rho(x,t)-0.7 \vert dx
$$
as  a measure for  the persistence of the periodic solution. One clearly observes that $I$ tends to zero in the  stable cases
considered here. The solution tends to the constant solution.
The convergence rate is larger for larger Chapman-Enskog term $D$, as expected.

\subsection{Test-case  3: Persistent periodic waves for  unstable situations}
This time we reconsider Example \ref{exN3}  using an equidistant spacing of the velocity points
$v_i = \frac{i}{N} $   and $N=20$. The relaxation parameter is again chosen as $\epsilon =0.01$.

Note that in this case choosing additionally  $\alpha_i = \frac{2 i}{N(N-1)}$ one  obtains

\begin{align*}
	\lambda =  \sum_{i=1}^{N-1} \alpha_i  v_i =  \sum_{i=1}^{N-1} \frac{2i^2}{N^2 (N-1)}
	= \frac{1}{3} \frac{2N-1}{N} 
\end{align*}	
and 
\begin{align*}
	\gamma =  \sum_{i=1}^{N-1} \frac{\alpha_i }{ v_i }=  \sum_{i=1}^{N-1} \frac{2i N}{i N (N-1)}
	= 2.
\end{align*}	
We obtain for the equilibrium functions 
\begin{align*}
	f_i^e(\rho,v) &= \frac{6}{N-1} \frac{N}{N+1} (F-E),\quad  i=1, \ldots, N-1\\
	f_N^e(\rho,v) &= F- 3 \frac{N}{N+1} (F-E)\\
	f_0^e(\rho,v) &= \rho- F-3 \frac{N-1}{N+1}(F-E).
\end{align*}	
To obtain an unstable situation we choose here the non-concave function $F= \rho(1-\rho)^2$  and $E= F(1-\frac{1}{3}\rho)$ and $E=F$, respectively.
The  functions $E$ fulfil the realizability condition  (\ref{realcont}).
We consider  the periodic   initial condition
$$
f_i (0,x) =  \frac{1}{N+1}(0.7+ 0.1\sin(6 \pi x)),\quad  i= 0,\ldots, N.
$$

\begin{figure}[h!]
		\tikzsetnextfilename{Test3_1}
		\begin{tikzpicture}[scale=0.7]
			\begin{axis}[xlabel =  $\rho$,
				legend style = {at={(0.5,1)},xshift=0.0cm,yshift=0.1cm,anchor=south},
				legend columns= 4,
				]
				\addplot[color = red,thick] file{datacont/Eu.txt};
				\addlegendentry{$E=F(1-\frac{\rho}{3})$}	
				\addplot[color = blue,thick] file{datacont/Fu.txt};
				\addlegendentry{$E=F$}
				\addplot[color = black, thick] file{datacont/Ru.txt};
				\addlegendentry{Real}
			\end{axis}
		\end{tikzpicture}
	\tikzsetnextfilename{Test3_2}
		\begin{tikzpicture}[scale=0.7]
			\begin{axis}[ylabel = $D$,xlabel =  $\rho$,
				legend style = {at={(0.5,1)},xshift=0.0cm,yshift=0.1cm,anchor=south},
				legend columns= 4,
				]
				
				\addplot[color = red,thick] file{datacont/unstabcond.txt};
				\addlegendentry{$E=F(1-\frac{\rho}{3})$}	
				\addplot[color = blue,thick] file{datacont/stabcond.txt};
				\addlegendentry{$E=F$}
				\addplot[black,thick] coordinates {(0,0)  (1,0)};
				{\tiny }
			\end{axis}
		\end{tikzpicture}
	
	\caption{Persistence of waves in unstable situations,  $N=20$. Left: functions $F= \rho(1-\rho)^2$  and $E= F(1-\frac{1}{3}\rho)$ and  realizability region. Right:   stability conditions for  $E= F(1-\frac{1}{3}\rho)$ and $E=F$.}
	\label{figunstable1}
\end{figure}
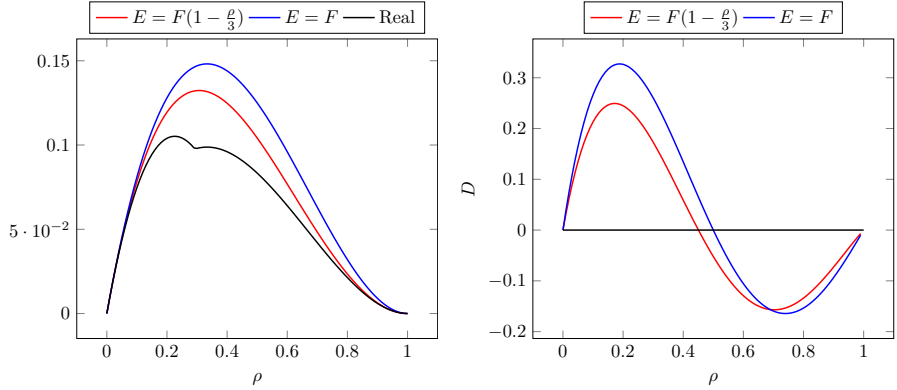

In Figure \ref{figunstable1} on the left, the functions $F$ and  $E=F(1-\frac{1}{3}\rho)$ and the lower bound of the  realizability domain are plotted.
On the right the Chapman -Enskog stability term $D(\rho)$ is plotted.

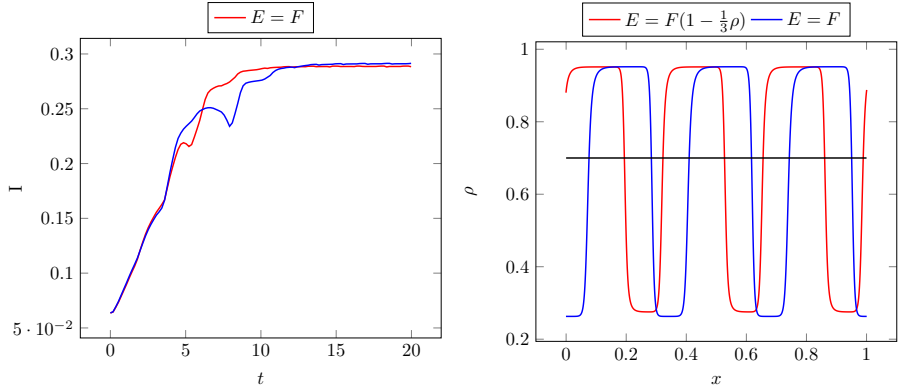
\begin{figure}[h!]
	\tikzsetnextfilename{Test3_3}
		\begin{tikzpicture}[scale=0.7]
			\begin{axis}[ylabel = I,xlabel =  $t$,
				legend style = {at={(0.5,1)},xshift=0.0cm,yshift=0.1cm,anchor=south},
				legend columns= 4,
				]
				
				\addplot[color = red,thick,each nth point={100}] file{data2/exp.txt};
				\addplot[color = blue,thick,each nth point={100}] file{data2/exp2.txt};
				\addlegendentry{$E=F$}
			\end{axis}
		\end{tikzpicture}
	\tikzsetnextfilename{Test3_4}
		\begin{tikzpicture}[scale=0.7]
			\begin{axis}[ylabel = $\rho$,xlabel =  $x$,
				legend style = {at={(0.5,1)},xshift=0.0cm,yshift=0.1cm,anchor=south},
				legend columns= 4,
				]
				\addplot[color = red,thick] file{data2/fig2.txt};
				\addlegendentry{$E= F(1-\frac{1}{3}\rho)$}
				\addplot[color = blue,thick] file{data2/fig22.txt};
				\addlegendentry{$E=F$}
				\addplot[black,thick] coordinates {(0,0.7)  (1,0.7)};
			\end{axis}
		\end{tikzpicture}
	
	\caption{Persistence of waves in  unstable situation with $F= \rho(1-\rho)^2$  and $E= F(1-\frac{1}{3}\rho)$, $E=F$, respectively. Left:  plot of $I$ versus time. Right:  plot of $\rho$ at final time $T=20$ for $N=20$.}
	\label{figunstable2}
\end{figure}

Figure \ref{figunstable2} shows on the left a plot of $I$ defined as above.
One observes that $I$ tends to a  positive value  for the  unstable cases considered in this subsection.
On the right,  the density $\rho$ is plotted at the final time $T=20$ showing the persistent periodic solutions for the two unstable cases considered here.

\subsection{Test-case 4:  Persistent periodic waves for a weakly unstable situation}

We consider again  Example \ref{exN2} with $N=2$.
As in the first example we consider the concave flux function 
$ F(\rho) = \rho (1-\rho)$. Here we  obtain  unstable situations by choosing appropriate $E$ leading to non-positive $D(\rho)$.
Thus, the relation beween flux and variance of the equilibrium distribution plays a major role in the appearance of persistent waves.

We consider two  numerical choices for the function $E$.  In both cases, $E$ is for smaller values of $\rho$ given by the flux function $F$. For larger values of  $\rho$ near $1$ we have chosen   $E_1$ to be  given by 
$E_1= 
\frac{1}{4} (\rho-1)(\rho-3)$ and $E_2$ by the lower boundary of the realizability domain.
In between a  $\mathcal{C}^2$-spline interpolation is used for both cases.

In Figure \ref{fig1b} on the left the functions $F,E_1$ and $E_2$  and the realizability domain  for $v_1 = \frac{1}{2}$ are plotted.
On the right the Chapman -Enskog stability terms are  plotted for the two choices of $E$.

The  initial conditions are again
$$
f_i (0,x) =  \frac{1}{3}(0.7+ 0.1\sin(6 \pi x)),\quad i= 0,1,2.
$$
Figure \ref{fig2} shows on the left a plot of 
$$
I (t)= \int_0^1 \vert \rho(x,t)-0.7 \vert dx
$$
as  a measure for  the persistence of the periodic solution. One observes that $I$ tends    positive values  for the unstable cases. 
On the right the density $\rho$ is plotted at the final time $T=20$ for the unstable situations.
Note that the  values og $I$ for large $t$ are smaller than the corresponding values in the previous  example.
The persistent waves are, correspondingly, less pronounced.
\begin{figure}[h!]
	\tikzsetnextfilename{Test4_1}
		\begin{tikzpicture}[scale=0.7]
			\begin{axis}[ylabel = $F / E$,xlabel =  $\rho$,
				legend style = {at={(0.5,1)},xshift=0.0cm,yshift=0.1cm,anchor=south},
				legend columns= 4,
				]
				\addplot[color = green,thick] file{data3/Fs.txt};
				\addlegendentry{$F$}
				\addplot[color = black, thick] file{data3/Rs.txt};
				\addlegendentry{Real}
				\addplot[color = blue,thick] file{data/E1.txt};
				\addlegendentry{$E1$}	
				\addplot[color = red,thick] file{data/E2.txt};
				\addlegendentry{$E2$}
			\end{axis}
		\end{tikzpicture}
	\tikzsetnextfilename{Test4_2}
		\begin{tikzpicture}[scale=0.7]
			\begin{axis}[ylabel = condition,xlabel =  $\rho$,
				legend style = {at={(0.5,1)},xshift=0.0cm,yshift=0.1cm,anchor=south},
				legend columns= 4,
				]
				
				\addplot[color = blue,thick] file{data/unstabcond1.txt};
				\addlegendentry{E1}	
				\addplot[color = red,thick] file{data/unstabcond2.txt};
				\addlegendentry{E2}
				\addplot[black,thick] coordinates {(0,0)  (1,0)};
			\end{axis}
		\end{tikzpicture}
	
	\caption{Persistence of waves in unstable situations, $N=2$. Left: functions $F$, $E_1$ and $E_2$ and  realizability region. Right:   stability conditions for  $E_1$ and $E_2$.}
	\label{fig1b}
\end{figure}
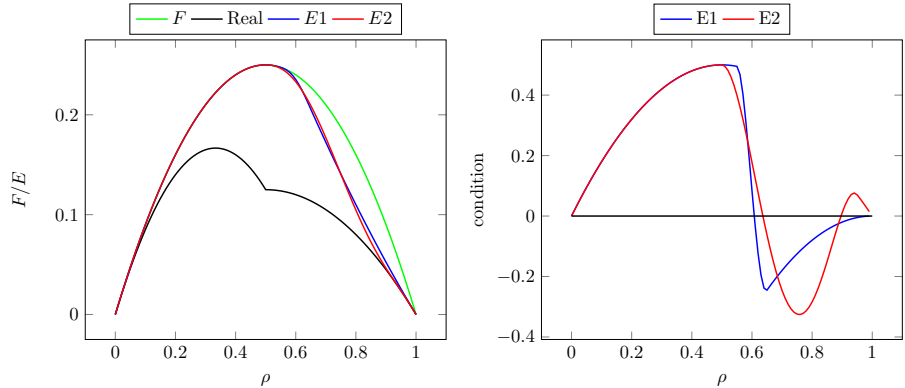

\begin{figure}[h!]
	\tikzsetnextfilename{Test4_3}
		\begin{tikzpicture}[scale=0.7]
			\begin{axis}[ylabel = I,xlabel =  $t$,
				legend style = {at={(0.5,1)},xshift=0.0cm,yshift=0.1cm,anchor=south},
				legend columns= 4,
				]
				\addplot[color = blue,thick] file{data/instable1a.txt};
				\addlegendentry{E1}
				\addplot[color = red,thick] file{data/instable1b.txt};
				\addlegendentry{E2}
			\end{axis}
		\end{tikzpicture}
	\tikzsetnextfilename{Test4_4}
		\begin{tikzpicture}[scale=0.7]
			\begin{axis}[ylabel = $\rho$,xlabel =  $x$,
				legend style = {at={(0.5,1)},xshift=0.0cm,yshift=0.1cm,anchor=south},
				legend columns= 4,
				]
				\addplot[color = blue,thick] file{data/fig1a.txt};
				\addlegendentry{E1}
				\addplot[color = red,thick] file{data/fig1b.txt};
				\addlegendentry{E2}
				\addplot[black,thick] coordinates {(0,0.7)  (1,0.7)};
			\end{axis}
		\end{tikzpicture}
	
	\caption{Persistence of waves in unstable situations. Left: Plot of $I$ versus time for $E_1$ and $E_2$. Right:   plot of the density $\rho$ for the  final time $T=20$ for $N=2$.}
	\label{fig2}
\end{figure}
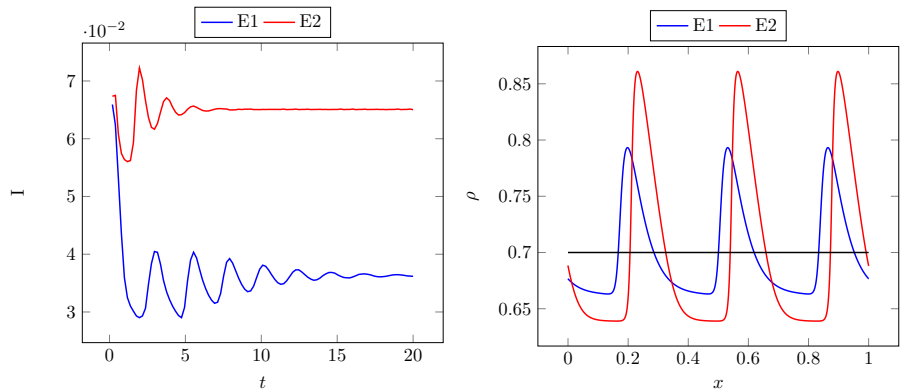

\section{Conclusions and outlook}

The paper presents a new class of non-linear kinetic discrete velocity models for traffic flow   having the correct invariant domain for traffic flow modelling and allowing for negative wave-speeds despite  the positivity of the microscopic velocities. The hyperbolic part is non-linear, but relatively simple,
being a totally linear degenerate and rich  hyperbolic problem with a simple structure of the integral curves.
For the homogeneous system, the fact, that the hyperbolic system is rich and totally linear degenerate, yields  existence of global entropy solutions,
see \cite{RS,Li}. For the relaxation system  the so-called semi-linear behaviour is obtained in \cite{CHN09} preventing the 
appearance of shock solutions.
Moreover, we have investigated conditions on the flux function $F$ such that  stop-and-go-wave like  instabilities arise, similar to  investigations that has been performed for the Aw-Rascle model in \cite{G04}.
Numerical results illustrate the behaviour of the solutions of Riemann problems wiout relaxation term and the persistence of periodic solutions for the full relaxation problem for various situations.
In particular, it turns out that the relation between flux  and variance of the equilibrium distribution plays a major role in the appearance of persistent waves.
Further research will be concerned with the development of higher order numerical methods for the kinetic  equations and a closer  investigation 
of the persistent periodic waves appearing for the relaxation problem.


\section{Appendix}
Consider  for  $k=0, \ldots, N-1$ the eigenvalues 
\begin{align*}
	\lambda_k (w) 
	&=v_k - \frac{ (v_{k+1}-v_k) w_{k+1}+ \sum_{j=k+2}^{N} (v_j-v_k) w_j \prod_{l=k+1}^{j-1}(1-w_l)}{\prod_{j=k+1}^{N}(1-w_j)}.
\end{align*}
For $j > k$ we prove (\ref{lamex}), i.e.
\begin{align*}
	\frac{\partial \lambda_k (w) }{\partial w_j } = - \frac{\lambda_j - \lambda_k}{1-w_j}.
\end{align*}
Consider first the case $j=k+1$, then
\begin{align*}
	\frac{\partial \lambda_k (w) }{\partial w_{k+1}} 
	&= - \frac{\partial }{\partial w_{k+1}} \left(\frac{ (v_{k+1}-v_k) w_{k+1}+ \sum_{j=k+2}^{N} (v_j-v_k) w_j \prod_{l=k+1}^{j-1}(1-w_l)}{\prod_{j=k+1}^{N}(1-w_j)} \right)\\
	&=  -\frac{ (v_{k+1}-v_k) \prod_{j=k+1}^{N}(1-w_j)+ (v_{k+1}-v_k) w_{k+1}  \prod_{j=k+2}^{N}(1-w_j)}{(\prod_{j=k+1}^{N}(1-w_j))^2}\\
	&=  -\frac{ (v_{k+1}-v_k) (1-w_{k+1})+ (v_{k+1}-v_k) w_{k+1}  }{(1-w_{k+1})^2\prod_{j=k+2}^{N}(1-w_j)}\\
	&=  -\frac{ v_{k+1}-v_k }{(1-w_{k+1})\prod_{j=k+1}^{N}(1-w_j)}\\
	&=- \frac{v_{k+1}-v_k }{1-w_{k+1}}  \frac{1-\sum_{j=0}^{k} f_j }{1-\rho} ,
\end{align*}
since
\begin{align*}
	&1-\sum_{j=0}^{k} f_j  =  1- w_0 -\sum_{j=1}^{k} w_j  \prod_{l=0}^{i-1} (1-w_l) \\
	&= (1-w_0) (1-w_1 -\sum_{j=2}^{k} w_j  \prod_{l=1}^{i-1} (1-w_l) )
	= \cdots = \prod_{j=0}^k (1-w_j).
\end{align*}
In the general case, we compute for $j>k+1$ 
\begin{align*}
	\frac{\partial \lambda_k (w) }{\partial w_{j}} 
	&= - \frac{\partial }{\partial w_{j}} \left(\frac{ (v_{k+1}-v_k) w_{k+1}+ \sum_{l=k+2}^{N} (v_l-v_k) w_l \prod_{r=k+1}^{l-1}(1-w_r)}{\prod_{l=k+1}^{N}(1-w_l)} \right)\\
	&=  -\frac{ (v_{j}-v_k) \prod_{r=k+1}^{j-1}(1-w_r)  - \frac{1}{1-w_j} \sum_{l=j+1}^N (v_l-v_k) w_l  \prod_{r=k+1}^{l-1}(1-w_r) }{\prod_{l=k+1}^{N}(1-w_l)} \\
	&- \frac{1}{1-w_j} \frac{ (v_{k+1}-v_k) w_{k+1}  +  \sum_{l=k+2}^{N}(v_l-v_k)w_l \prod_{r=k+1}^{l-1}(1-w_r)}{\prod_{l=k+1}^{N}(1-w_l)}
\end{align*}
This yields 
\begin{align*}
	\frac{\partial \lambda_k (w) }{\partial w_{j}} 
	&=  -\frac{ (v_{j}-v_k) \prod_{r=k+1}^{j-1}(1-w_r)   }{\prod_{l=k+1}^{N}(1-w_l)} - \frac{1}{1-w_j} \frac{ (v_{k+1}-v_k) w_{k+1}  }{\prod_{l=k+1}^{N}(1-w_l)}\\ 
	&+
	\frac{ \frac{1}{1-w_j} \sum_{l=j+1}^N (v_l-v_k) w_l  \prod_{r=k+1}^{l-1}(1-w_r) }{\prod_{l=k+1}^{N}(1-w_l)} \\
	&- \frac{1}{1-w_j} \frac{   \sum_{l=k+2}^{N}(v_l-v_k)w_l \prod_{r=k+1}^{l-1}(1-w_r)}{\prod_{l=k+1}^{N}(1-w_l)}\\
	&=  -\frac{ (v_{j}-v_k) \prod_{r=k+1}^{j-1}(1-w_r)   }{\prod_{l=k+1}^{N}(1-w_l)} - \frac{1}{1-w_j} \frac{ (v_{k+1}-v_k) w_{k+1}  }{\prod_{l=k+1}^{N}(1-w_l)}\\ 
	&-\frac{1}{1-w_j} 
	\frac{ \sum_{l=k+2}^j (v_l-v_k) w_l  \prod_{r=k+1}^{l-1}(1-w_r) }{\prod_{l=k+1}^{N}(1-w_l)}
\end{align*}
and then
\begin{align*}
	\frac{\partial \lambda_k (w) }{\partial w_{j}} 
	&=  -\frac{1}{1-w_j}  \frac{ (v_{j}-v_k) \prod_{r=k+1}^{j}(1-w_r)   }{\prod_{l=k+1}^{N}(1-w_l)} - \frac{1}{1-w_j} \frac{ (v_{k+1}-v_k) w_{k+1}  }{\prod_{l=k+1}^{N}(1-w_l)}\\ 
	&-\frac{1}{1-w_j} 
	\frac{ \sum_{l=k+2}^j (v_l-v_k) w_l  \prod_{r=k+1}^{l-1}(1-w_r) }{\prod_{l=k+1}^{N}(1-w_l)} \\
	&=- \frac{1}{1-w_{j}} \left((v_{j}- v_k)  \frac{1}{\prod_{r=j+1}^N (1-w_r)}+ \sum_{l=k+1}^{j} (v_l -v_k)w_l \frac{1}{\prod_{r=l}^{N} (1-w_r)}\right)\\
	&=- \frac{1}{1-w_{j}} \left((v_{j}- v_k)  \frac{\prod_{l=0}^{j} (1-w_l)}{\prod_{r=0}^N (1-w_r)}- \frac{\sum_{l=k+1}^{j} (v_l -v_k)w_l \prod_{r=0}^{l-1} (1-w_r)}{\prod_{r=0}^{N} (1-w_r)}\right)\\
	&=- \frac{1}{1-w_{j}} \left((v_{j}- v_k)  \frac{1-\sum_{l=0}^{j} f_l}{1-\rho}+ \frac{\sum_{l=k+1}^{j} (v_l -v_k)f_l}{1-\rho}\right).
\end{align*}
On the other hand we have for $j>k$
\begin{align*}
	&\lambda_{j} (w) - \lambda_k (w) \\&= v_{j}- v_k- \frac{1}{1-\rho}\left(
	\sum_{l=j+1}^{N} (v_l-v_{j}) f_l-  \sum_{l=k+1}^{N} (v_l-v_{k}) f_l\right)\\
	&=v_{j}- v_k- \frac{1}{1-\rho}\left(
	\sum_{l=j+1}^{N} (v_k -v_j) f_l- \sum_{l=k+1}^{j}  (v_{l}-v_{k}) f_{l} \right)\\
	&=(v_{j}- v_k) \left(1+ \frac{1}{1-\rho} \sum_{l=j+1}^{N} f_l \right)
	+\frac{1}{1-\rho}\sum_{l=k+1}^{j}  (v_{l}-v_{k}) f_{l} \\
	&=(v_{j}- v_k)  \frac{1-\sum_{l=0}^{j} f_l}{1-\rho} +\frac{1}{1-\rho}\sum_{l=k+1}^{j}  (v_{l}-v_{k}) f_{l}  \\
\end{align*}
Altogether we have proven (\ref{lamex}) for $j>k$.

\end{document}